\documentclass[twocolumn,10pt]{asme2e_m}

\usepackage[cmex10]{amsmath}
\usepackage{amsmath,amssymb,amsfonts}
\usepackage{amsbsy}
\usepackage[table]{xcolor}
\usepackage{tikz}
\usepackage{gensymb}
\usepackage{latexsym}
\usepackage{soul}
\usepackage{array}
\usepackage{cleveref}
\usepackage{color}
\usepackage{relsize}
\usepackage{tikz}
\usepackage{cite}
\usepackage{epsfig}
\usepackage{fancyhdr}
\usepackage{graphics}
\usepackage{float}
\usepackage{tablefootnote}
\usepackage{breqn} 
\usepackage{eucal}

\confshortname{ASME 2017 Dynamic Systems and Control Conference} \conffullname{DSCC2017} \confdate{11-13}
\confmonth{October} \confyear{2017} \confcity{Tysons, Virginia} \confcountry{USA}
\papernum{DSCC2017-5010}
\title{MIMO First and Second Order Discrete Sliding Mode Controls of Uncertain Linear Systems under Implementation Imprecisions}
\author{Mohammad Reza Amini\thanks{Address all correspondence to this author.} 
       \affiliation{
Dept. of Mechanical Engineering\\
Michigan Technological University\\
Houghton, MI 49931\\
Email: mamini@mtu.edu
    }
}
\author{Mahdi Shahbakhti
       \affiliation{
Dept. of Mechanical Engineering\\
Michigan Technological University\\
Houghton, MI 49931\\
Email: mahdish@mtu.edu
    }
}
\author{Selina Pan
    \affiliation{Research and Innovation Center\\
Ford Motor Company\\
Palo Alto, CA, 94304\\
Email: span6@ford.com
    }
}
\graphicspath{{./figures/}}
\begin{document}
\pagestyle{empty}
\maketitle
\thispagestyle{empty}
\pagestyle{empty}
\vspace{-1cm}
\begin{abstract}
The performance of a conventional model-based controller significantly depends on the accuracy of the modeled dynamics. The model of a plant's dynamics is subjected to errors in estimating the numerical values of the physical parameters, and variations over operating environment conditions and time. These errors and variations in the parameters of a model are the major sources of uncertainty within the controller structure. Digital implementation of controller software on an actual electronic control unit (ECU) introduces another layer of uncertainty at the controller inputs/outputs. The implementation uncertainties are mostly due to data sampling and quantization via the analog-to-digital conversion (ADC) unit. The failure to address the model and ADC uncertainties during the early stages of a controller design cycle results in a costly and time consuming verification and validation (V\&V) process. In this paper, new formulations of the first and second order discrete sliding mode controllers (DSMC) are presented for a general class of uncertain linear systems. The knowledge of the ADC imprecisions is incorporated into the proposed DSMCs via an online ADC uncertainty prediction mechanism to improve the controller robustness characteristics. Moreover, the DSMCs are equipped with adaptation laws to remove two different types of modeling uncertainties (multiplicative and additive) from the parameters of the linear system model. The proposed adaptive DSMCs are evaluated on a DC motor speed control problem in real-time using a processor-in-the-loop (PIL) setup with an actual ECU. The results show that the proposed SISO and MIMO second order DSMCs improve the conventional SISO first order DSMC tracking performance by 69\% and 84\%, respectively. Moreover, the proposed adaptation mechanism is able to remove the uncertainties in the model by up to 90\%.
\end{abstract}
\vspace{-0.5cm}
\begin{nomenclature}\hspace{-8.5mm} 
\begin{tabular*}{0.95\columnwidth}{ l p{8cm}}
\multicolumn{2}{l}{   }  \vspace{-2.5mm} \\ 
$x,~\mathbf{X}$ &states of the system [$-$]\\ 
$s,~\mathbf{S}$ &first order sliding surface [$-$]\\ 
$\xi,~\boldsymbol{\Xi}$ &second order sliding surface [$-$]\\
$J$ &effective DC motor's rotor inertia [$kg.m^2$]\\ 
\end{tabular*}
\end{nomenclature}
%
\begin{table} \vspace{-2.5mm}
{\renewcommand{\arraystretch}{1.0}
\begin{tabular*}{0.01 \columnwidth}{ p{1cm} p{8cm}}
\multicolumn{2}{l}{   } \\
$R$ & electrical resistance, [$\Omega$]\\
$L$ & electrical inductance, [H]\\
$\Gamma$ & generated torque on DC motor's rotor, [N.m]\\
$k_f$ & mechanical damping, [N.m.s]\\
$k_m$ & motor torque constant, [N.m/A]\\
$k_b$ & electromotive force constant, [V.s/rad]\\
$T$ &sampling time, [$s$]\\ 
$\rho,~\mathcal{P}$ & tunable first order DSMC gain/matrix, [-]\\
$\varphi,~\boldsymbol{\Phi}$ & tunable second order DSMC gain/matrix, [-]\\
$\rho_{\beta}$ &adaptation gain for multiplicative uncertainty, [$-$]\\ 
$\rho_{\alpha}$ &adaptation gain for additive uncertainty, [$-$]\\ 
$\alpha$ &additive uncertainty term, [$-$]\\ 
$\beta$ &multiplicative uncertainty term, [$-$]\\ 
\end{tabular*}
} \vspace{-0.5cm}
\end{table}

\section{INTRODUCTION} \label{Sec:Intro}
There are two major sources of uncertainties which make the completion of traditional verification and validation (V$\&$V) cycle of a model-based controller challenging and costly. These uncertainties are mostly due to: (i) errors in estimating the model parameters and physical changes in the plant or fluctuations in the environment in which the system operates, and (ii) imprecisions often arise during digital implementation of the controller software on an actual electronic control unit (ECU) via analog-to-digital converter (ADC) unit. Neglecting the modeling and implementation uncertainties during the early stages of the controller design cycle leads to substantial deviation in the controller performance once it is implemented in the real ECU~\cite{NASA_Dabney,shahbakhti2012early}. 

It has been shown in the literature~\cite{Amini_CEP,Amini_DSCC2016,Pan_DSC} that the performance of a conventional model-based controller is considerably sensitive to any errors in the modeled plant dynamics, and ignoring the uncertainties in the model results in the controller failure. Moreover, digital implementation of a controller software on an actual ECU introduces sampling and quantization imprecision on the controller input/output (I/O) signals via the ADC unit. The impact of implementation imprecisions on a conventional model-based controller has been studied in the literature~\cite{AminiSAE2016,Hansen_DSCC}, and it has been shown that the ADC imprecisions can make a controller significantly deviates from its desired response. 

One effective approach to minimize the gap between the designed and implemented controllers is early model-based design and verification of the controller software~\cite{Amini_DSC}. In this approach, prior to conducting the conventional V$\&$V iterative cycle, the structure of the model-based controller is investigated to identify the uncertainty sources. Next, the controller structure is modified to achieve higher robustness, and have adaptability against uncertainties within the model and implementation imprecisions. Among different model-based controller design techniques, sliding mode control (SMC)~\cite{AminiSAE2016} and discrete sliding mode control (DSMC)~\cite{Pan_DSC,Kyle_CDC,Amini_DSC} have shown to be low-cost solutions, in which their structures allow for achieving higher robustness characteristics against model and ADC uncertainties. 

\begin{figure}[h!]
\begin{center}
\includegraphics[angle=0,width= \columnwidth]{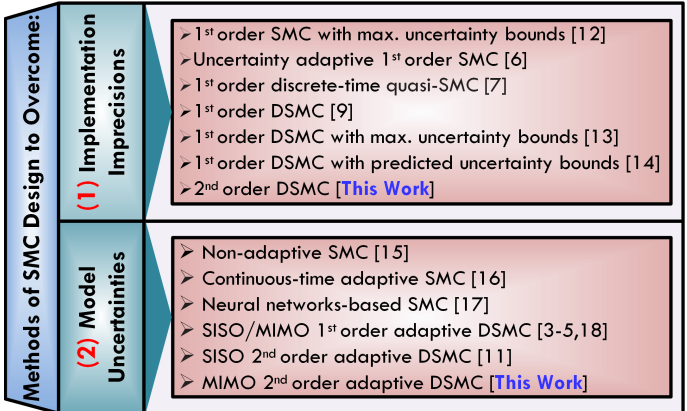} \vspace{-0.75cm}
\caption{\label{fig:Survay_Sampling_ACC2017}Background of previous SMC studies\cite{KyleACC,AminiSAE2016,Hansen_DSCC,Kyle_CDC,Kyle_DSCC,Amini_ACC2016,Misawa_DSC,Slotine,fang1999use,Chan_Automatica,Pan_DSC,Amini_DSCC2016,Amini_CEP,Amini_ACC2017} on controller design
against ADC and modeling uncertainties.} \vspace{-1.0cm}
\end{center}
\end{figure}

Figure~\ref{fig:Survay_Sampling_ACC2017} shows the previous works in the literature which aimed to provide a solution to handle modeling and ADC uncertainties based on the SMC/DSMC structure. The concept of the second order DSMC~\cite{mihoub2009real} helps to reduce the high frequency oscillation due to chattering, and improves the controller robustness against data sampling imprecisions~\cite{Amini_ACC2017}. This is because in the second order DSMC, not only is the system state driven to the sliding manifold, but the state derivative (difference function) is also steered to zero. In this paper, we employ an adaptive robust DSMC design from~\cite{Amini_CEP,Amini_ACC2017} to remove the uncertainties in the model and improve the robustness characteristics against ADC imprecisions. The generic single-input single-output (SISO) adaptive DSMC formulation from~\cite{Amini_CEP,Amini_ACC2017} is extended to a general class of uncertain linear systems, and will be used to formulate SISO and multi-input multi-output (MIMO) first and second order adaptive DSMCs with predicted ADC uncertainties. This paper presents the first theoretical development for application of the \textit{MIMO second order DSMC for uncertain linear systems} under implementation imprecisions. \vspace{-0.7cm}

The contribution of this paper is twofold. First, new first and second order DSMCs with incorporated implementation imprecisions are formulated for a general class of SISO/MIMO linear systems under ADC uncertainties. Second, the proposed controller design is extended to handle model uncertainties using a discrete Lyapunov stability argument that also guarantees the asymptotic stability of the closed-loop system. The application of the proposed first and second order adaptive DSMCs are shown on an uncertain model of a DC motor for speed control.\vspace{-0.65cm}
%
\section{Adaptive Robust DSMC for Linear Systems} \label{sec:MainTheory}
In the absence of model uncertainties and implementation imprecisions, the state-space equation of a $r^{~th}$-order linear system, discretized by using a first order Euler approximation~\cite{Hansen_DSCC}, is represented by:\vspace{-0.4cm}
\begin{equation}\label{eq:Generic_linear_1}
\mathbf{X}_{r\times 1}(i+1)= \left(T\mathbf{A}_{r\times r}+\mathbf{I}_{r\times r}\right)\mathbf{X}_{r \times 1}(i)+T\mathbf{B}_{r\times h}\mathbf{U}_{h \times 1}\vspace{-0.25cm}
\end{equation} 
where $x\subset{\mathbf{X}\in{\mathbb{R}^{r}}}$, $u\subset{\mathbf{U}\in{\mathbb{R}^{h}}}$, and ${\mathbf{I}\in{\mathbb{R}^{r\times r}}}$ are the state vector, the control input vector, and identity matrix, respectively. Matrix $\mathbf{A}$ represents the linear system dynamics, and its elements, $a_{pq},~p,q=1,...,r$ are obtained based on the system's physical equations and the interaction among different states. In practice, {assuming an exact knowledge of $\mathbf{B}$}, the identified matrix $\mathbf{A}$ is subjected to several sources of uncertainties, e.g. the plant parameter variations over time. 

The linear system in Eq.~(\ref{eq:Generic_linear_1}), in the presence of additive ($\alpha$) and multiplicative ($\beta$) types of modeling uncertainties, and ADC imprecisions, can be presented as follows: \vspace{-0.4cm}
\begin{equation}\label{eq:Generic_linear_2}
\begin{split}
\mathbf{X}_{r \times 1}(i+1)=~~~~~~~~~~~~~~~~~~~~~~~~~~~~~~~~~~~~~~~~~~~~~~~~~~~~~~~~~~~~~~~~~~~~~\\
\left(T\begin{bmatrix}
\beta_{11}a_{11}+\alpha_{11}~~ & \beta_{12}a_{12}+\alpha_{12} & \cdots & \beta_{1r}a_{1r}+\alpha_{1r} \\ 
\beta_{21}a_{21}+\alpha_{21}~~ & \beta_{22}a_{22}+\alpha_{22} & \cdots & \beta_{2r}a_{2r}+\alpha_{2r}\\ 
\vdots & \vdots & \ddots & \vdots \\ 
\beta_{r1}a_{r1}+\alpha_{r1}~~ & \beta_{r2}a_{r2}+\alpha_{r2} & \cdots & \beta_{rr}a_{rr}+\alpha_{rr} 
\end{bmatrix} \right.\\
\left.+\begin{bmatrix}
1 & 0 & \cdots & 0 \\ 
0 & 1 & \cdots & 0\\ 
\vdots & \vdots & \ddots & \vdots \\ 
0 & 0 & \cdots & 1
\end{bmatrix}_{r \times r}\right)
\begin{bmatrix}
x_1(i)\\ 
x_2(i)\\ 
\vdots\\ 
x_r(i)
\end{bmatrix}+T\mathbf{B}_{r \times h}
\left(\mathbf{U}_{h \times 1}+\delta_{h \times 1}(i)
\right)	\vspace{-0.35cm}
\end{split}
\end{equation} 
where $a_{pq}$ represents nominal values of the ideal linear system's dynamics ($\mathbf{A}$), and $\delta_{h \times 1}$ is the vector of propagated ADC uncertainties on the control signal ($\mathbf{U}$). As can be seen, for each element of $\mathbf{A}$, one additive ($\alpha$) and one multiplicative ($\beta$) unknown terms are considered to present any errors or variation in the values of the model parameters. For the uncertain linear system in Eq.~(\ref{eq:Generic_linear_2}) under implementation imprecisions, the tracking control problem is defined to drive the states of the system ($x$) to their desired values ($x_d\subset{\mathbf{X}_d\in{\mathbb{R}^{r}}}$). To this end, SISO/MIMO first and second order DSMCs, under modelling and implementation uncertainties, are formulated in the following sections. \vspace{-0.5cm}

\subsection{First Order Adaptive DSMC} \label{FirstOrderDSMC}
For the linear system in Eq.~(\ref{eq:Generic_linear_1}), the vector ($\mathbf{S}$) of the first order sliding surface variables ($s$) is defined as the difference between the desired ($x_d$) and the measured signal ($x$) as follows:
\vspace{-0.65cm}
\begin{gather}\label{eq:StageIII_9}
\mathbf{S}(i)=\mathbf{X}(i)-\mathbf{X}_d(i)
\vspace{-0.5cm}
\end{gather}
%
The objective is asymptotic and finite-time convergence of $s$ to zero. To achieve this, the vector of the first order DSMC input $\mathbf{U}(i)$ is obtained according to the following first order sliding reaching law~\cite{Kyle_CDC,Amini_ACC2016}:
\vspace{-0.5cm}
\begin{gather}\label{eq:StageIII_10}
|\mathbf{S}(i+1)|\leq \mathbf{\mathcal{P}} |\mathbf{S}(i)| \vspace{-0.5cm}
\end{gather}
in which, $\mathbf{\mathcal{P}}$ is the matrix of tunable DSMC gains. For a SISO DSMC, $\mathbf{\mathcal{P}}$ is diagonal:~$\mathbf{\mathcal{P}}=diag[\rho_{1},...,\rho_r]$, where $0<\rho_{1,...,r}<1$~\cite{Kyle_CDC}. On the other hand, for a MIMO DSMC structure, the off-diagonal elements of $\mathbf{\mathcal{P}}$ can be non-zero; however, the eigenvalues of $\mathbf{\mathcal{P}}$ should lie within the unit circle to guarantee the closed-loop system stability~\cite{Pan_DSC}. If $r=h$, $\mathbf{B}$ is a square matrix, and the control input vector $\mathbf{U}$ can be calculated as follows for the linear system according to the sliding reaching law from Eq.~(\ref{eq:StageIII_10}), {assuming that $\mathbf{B}$ is invertible}:
\vspace{-0.4cm}
\begin{gather}\label{eq:StageIII_12_1}
\mathbf{U}(i)=
\mathbf{B}^{-1}\left(\frac{1}{T}\left[({\mathbf{\mathcal{P}}}-\mathbf{I})\mathbf{X}(i)-\mathbf{\mathcal{P}} \mathbf{X}_d(i)+\mathbf{X}_d(i+1)\right]-\mathbf{A}\mathbf{X}(i)\right) \vspace{-0.5cm}
\end{gather}

In the absence of model uncertainties and implementation imprecisions, Eq.~(\ref{eq:StageIII_12_1}) calculates the control input vector of a SISO/MIMO first order DSMC at each time step, in which for each state variable, a sliding surface is defined and it is assumed that a unique control input, either physical or synthetic, exists for every single sliding surface. As discussed earlier in Section~\ref{sec:MainTheory}, the signals at the controller I/O are subjected to sampling and quantization imprecisions. The introduced ADC imprecisions on the measured signals ($\mu_x$) are next propagated through the state space equations of the linear system. The propagated ADC imprecisions on the control signals are shown in Eq.~(\ref{eq:Generic_linear_2}) by $\delta$. Here, the approach to overcome the ADC uncertainties is to estimate the overall uncertainty bounds on the control inputs in real-time, and then, modify the DSMC to make it more robust against sampling and quantization imprecisions with respect to the knowledge of the propagated ADC uncertainties.

For the linear system under ADC imprecisions, unlike nonlinear systems~\cite{Amini_ACC2016}, the propagated uncertainties on the control signal ($\delta$) can be calculated analytically with respect to the uncertainties on the measured signals ($\mu_x$). This means that if the introduced imprecisions on the measured signals ($\mu_x$) could be predicted online, $\delta$ can be estimated in real-time. Here, it is assumed that the state and control input of an ideal DSMC, where there are neither model uncertainties nor ADC imprecisions, are shown by $\bar{\mathbf{X}}$ and $\bar{\mathbf{U}}$, respectively. According to Eq.~(\ref{eq:StageIII_12_1}), $\bar{\mathbf{U}}$, the control actions of the ideal DSMC, can be established as:\vspace{-0.5cm}
\begin{gather}\label{eq:StageIII_13}
\bar{\mathbf{U}}(i)=
\mathbf{B}^{-1}\left(\frac{1}{T}[(\mathcal{P}-\mathbf{I})\bar{\mathbf{X}}(i)-\mathcal{P} \mathbf{X}_d(i)+\mathbf{X}_d(i+1)]-\mathbf{A}\bar{\mathbf{X}}(i)\right) \vspace{-0.5cm}
\end{gather}

The only difference between the ideal DSMC (Eq.~(\ref{eq:StageIII_13})) and DSMC under ADC imprecisions (Eq.~(\ref{eq:StageIII_12_1})) is the introduced sampling and quantization uncertainties. Thereby, the difference in the control input vectors from Eq.~(\ref{eq:StageIII_12_1}) and~(\ref{eq:StageIII_13}) is assumed to be the propagated uncertainty on control input. The difference between $\mathbf{U}$ and $\bar{\mathbf{U}}$ can be analytically found with respect to Eq.~(\ref{eq:StageIII_12_1}) and~(\ref{eq:StageIII_13}): \vspace{-0.5cm}
%
%
\begin{gather}\label{eq:StageIII_15_1}
\bar{\mathbf{U}}(i)-{\mathbf{U}}(i)=\mathbf{B}^{-1}\left(\frac{1}{T}[(\mathcal{P}-\mathbf{I})\mu_{_X}(i)]-\mathbf{A}\mu_{_X}(i)\right)\vspace{-0.5cm}
\end{gather}
where, $\mu_{_X}=\bar{\mathbf{X}}-\mathbf{X}=\{\mu_{x_1},\mu_{x_2},...,\mu_{x_r}\}$ is the introduced ADC imprecisions on the measured signals. In practice, $\mu_{_X}$ cannot be always obtained, because there is not usually access to the measured signals before the ADC unit, and the only available information is the discretized and digitized signal after the ADC. In our previous work~\cite{Amini_ACC2016}, we proposed a simple and accurate technique to estimate $\mu_{_X}$ with respect to the measured signal after ADC ($\mathbf{X}$), the change rate in the measured signal (slope of the signal), sampling time ($T$), and the quantization level ($n$). 

Overall, the estimated sampling and quantization uncertainties ($\hat{\mu}_{x}$) on a measured signal at $i^{th}$ time step can be obtained by using the following equation:\vspace{-0.5cm}
\begin{gather}\label{eq:StageIII_7}
\hat{\mu}_{x}(i)=x(i)-x(i-1)+\frac{1}{2} \frac{FSR}{2^n}
\end{gather}
where, $FSR$ is the full scale range of the measured signal, $n$ is number of ADC bits and represents the ADC resolution, and $\hat{\mu}_{x}$ is the estimated ADC uncertainty on the measured signal, and it is assumed to be the difference between the measured signals before and after ADC ($\hat{\mu}_{x}=\bar{x}-x$), where $\bar{x}$ is the estimated actual measured signal before ADC. Upon substituting the $\hat{\mu}_{x}$ into Eq.~(\ref{eq:StageIII_15_1}) at each time step, the propagated ADC uncertainties on the control signals can be obtained.

A new diagonal matrix ($\mu_{\mathbf{U}}$) is defined with respect to the estimated propagated ADC uncertainties on the control signals:\vspace{-0.4cm}
\begin{gather}\label{eq:StageIII_14}
\mu_{\mathbf{U}}=diag[\bar{u}_1-u_1,...,\bar{u}_r-u_r]
\end{gather}
where $u\subset{\mathbf{U}}$ and $\bar{u}\subset{\bar{\mathbf{U}}}$. The diagonal elements of $\mu_{\mathbf{U}}$ are calculated with respect to Eq.~(\ref{eq:StageIII_15_1}). According to~\cite{Amini_ACC2016,Amini_CEP,Amini_DSCC2016}, by inclusion of the propagated ADC uncertainties on control signals ($\mu_{\mathbf{U}}$) into the DSMC structure, the robustness of the conventional sliding mode controller can be improved. For the linear system (Eq.~(\ref{eq:StageIII_12_1})), in the absence of model uncertainties, the control input ($^{1DSMC}\mathbf{U}$) of the conventional first order DSMC is modified against data sampling and quantization imprecisions as follows:\vspace{-0.4cm}
\begin{gather}\label{eq:StageIII_21}
^{1DSMC}\mathbf{U}^{mod}(i)=~~~~~~~~~~~~~~~~~~~~~~~~~~~~~~~~~~~~~~~~~~~~~~~~~~~~~~~~~~~~~~~\\
\mathbf{B}^{-1}\Bigg(\frac{1}{T}\{(\mathcal{P}-\mathbf{I}){\mathbf{X}}(i)-\mathcal{P} \mathbf{X}_d(i)+\mathbf{X}_d(i+1)\}-\mathbf{A}\mathbf{X}(i)\Bigg)\nonumber\\
-|\mu_{\mathbf{U}}(i)|\times sat(\mathbf{S}(i)) \nonumber
\end{gather}
where $\mu_{\mathbf{U}}$ is calculated according to Eq.~(\ref{eq:StageIII_14}), and $sat(.)$ is the saturation function which is used instead of the $signum$ function to avoid possible high frequency chattering which occurs
in discrete systems during implementation of the $signum$
function~\cite{Slotine,Amini_DSC}.

In the next step, the uncertainties in the model are included in the first order DSMC formulation. It can be easily shown that the first order sliding vector ($\mathbf{S}$) for the linear system with unknown additive and multiplicative parameters becomes: \vspace{-0.4cm}
\begin{equation}\label{eq:Generic_linear_3}
\begin{split}
\mathbf{S}_{r\times 1}(i+1)=\mathcal{P}_{r\times r} \mathbf{S}_{r\times 1}(i)+~~~~~~~~~~~~~~~~~~~~~~~~~~~~~~~~~~~~~~~~~~~\\
T\begin{bmatrix}
\tilde{\beta}_{11} a_{11}+\tilde{\alpha}_{11}~~ & \tilde{\beta}_{12} a_{12}+\tilde{\alpha}_{12} & \cdots & \tilde{\beta}_{1r} a_{1r}+\tilde{\alpha}_{1r} \\ 
\tilde{\beta}_{21} a_{21}+\tilde{\alpha}_{21}~~ & \tilde{\beta}_{22} a_{22}+\tilde{\alpha}_{22} & \cdots & \tilde{\beta}_{2r} a_{2r}+\tilde{\alpha}_{2r}\\ 
\vdots & \vdots & \ddots & \vdots \\ 
\tilde{\beta}_{r1} a_{r1}+\tilde{\alpha}_{r1}~~ & \tilde{\beta}_{r2} a_{r2}+\tilde{\alpha}_{r2} & \cdots & \tilde{\beta}_{rr} a_{rr}+\tilde{\alpha}_{rr} 
\end{bmatrix}\mathbf{X}_{r\times 1}(i)
\end{split}
\end{equation} 
where $\tilde{\beta}_{pq}=\beta_{pq}-\hat{\beta}_{pq}$ and $\tilde{\alpha}_{pq}=\alpha_{pq}-\hat{\alpha}_{pq}$ are the errors in estimating the unknown \textit{constant} multiplicative and additive parameters, respectively. 

$\blacksquare$~\textbf{Theorem}: The adaptation laws for converging the unknown parameters of the uncertain linear system (Eq.~(\ref{eq:Generic_linear_2})) to their nominal values ($\hat{\beta}\rightarrow \beta$,~$\hat{\alpha}\rightarrow \alpha$), based on the DSMC formulation, are as follows:\vspace{-0.6cm}
\begin{subequations} \label{eq:Generic_linear_4}
\begin{align}
\hat{\beta}_{pq}(i+1)=\hat{\beta}_{pq}(i)+\frac{Ts_p(i)a_{pq}x_q(i)}{\rho_{\beta_{pq}}} \\
\hat{\alpha}_{pq}(i+1)=\hat{\alpha}_{pq}(i)+\frac{Ts_p(i)x_q(i)}{\rho_{\alpha_{pq}}} 
\end{align} 
\end{subequations} 
where, $p,q=1...r$, and $\rho_{\beta}$ and $\rho_{\alpha}$ are tunable \textit{positive} multiplicative and additive adaptation gains chosen for the numerical sensitivity in the estimations of the unknown parameters. 

\textbf{Proof}: A Lyapunov stability analysis is performed to derive the adaptation laws and guarantee the stability of the closed-loop system. To this end, the analysis begins with a first-order linear uncertain system ($r=1$). The following positive definite scalar Lyapunov function ($\mathbf{V}$) for a first order linear system is proposed: \vspace{-0.6cm}
\begin{gather} \label{eq:DSCC2017_1}
V(i)=\frac{1}{2}s^2(i)+\frac{1}{2}\rho_{\beta}\tilde{\beta}^2(i)+\frac{1}{2}\rho_{\alpha}\tilde{\alpha}^2(i)
\end{gather}
As can be seen, $\mathbf{V}$ is a quadratic function of the tracking error ($\mathbf{S}$), and unknown parameters estimations. In the next step, the first order Lyapunov difference function ($\Delta V$) is calculated. $\Delta V$ is obtained by applying a Taylor series expansion on Eq.~(\ref{eq:DSCC2017_1}): \vspace{-0.4cm}
\begin{gather}\label{eq:DSCC2017_2}
\Delta V(i)=\frac{\partial V(i)}{\partial s(i)}\Delta s(i)+\frac{\partial V(i)}{\partial \tilde{\beta}(i)}\Delta \tilde{\beta}(i)+\frac{\partial V(i)}{\partial \tilde{\alpha}(i)}\Delta \tilde{\alpha}(i)+\\ 
\frac{1}{2} \frac{\partial^2 V(i)}{\partial {s}^2(i)}\Delta {s}^2(i)+\frac{1}{2} \frac{\partial^2 V(i)}{\partial {\tilde{\beta}}^2(i)}\Delta {\tilde{\beta}}^2(i)+
\frac{1}{2} \frac{\partial^2 V(i)}{\partial {\tilde{\alpha}}^2(i)}\Delta {\tilde{\alpha}}^2(i)+... \nonumber
\end{gather}
where, $\Delta s(i)\equiv s(i+1)-s(i),~\Delta \tilde{\beta}(i) \equiv \tilde{\beta}(i+1)-\tilde{\beta}(i),~\Delta \tilde{\alpha}(i) \equiv \tilde{\alpha}(i+1)-\tilde{\alpha}(i)$. Upon substitution of partial derivatives into Eq.~(\ref{eq:DSCC2017_2}), we have: \vspace{-0.5cm}
\begin{gather}\label{eq:DSCC2017_3}
\Delta V(i) 
=s(i)\Delta s(i)+ \rho_{\beta}\tilde{\beta}(i)\Delta \tilde{\beta}(i)+\rho_{\alpha}\tilde{\alpha}(i)\Delta \tilde{\alpha}(i)+\\
\frac{1}{2}\Delta {s}^2(i)+\frac{1}{2}\rho_{\beta}\Delta {\tilde{\beta}}^2(i)+\frac{1}{2}\rho_{\alpha}\Delta {\tilde{\alpha}}^2(i)+... \nonumber
\end{gather}
It is assumed that for small enough sampling periods, $\Delta {s}^2(i),~\Delta\tilde{\beta}^2(i),~\Delta\tilde{\alpha}^2(i)\approx 0$~\cite{Amini_DSC}. The same assumption is valid for higher order terms also ($>2$). By using Eq.~(\ref{eq:Generic_linear_3}) for the first order system, Eq.~(\ref{eq:DSCC2017_3}) can be simplified as follows:\vspace{-0.5cm}
\begin{gather}\label{eq:DSCC2017_4}
\Delta V(i) 
=s(i)\Big((\rho-1)s(i)+T(\tilde{\beta}(i)a+\tilde{\alpha}(i))x(i)\Big)+\\
\rho_{\beta}\tilde{\beta}(i)\Delta \tilde{\beta}(i)+\rho_{\alpha}\tilde{\alpha}(i)\Delta \tilde{\alpha}(i) \nonumber
\end{gather}

Eq.~(\ref{eq:DSCC2017_4}) can be re-arranged as follows:\vspace{-0.5cm}
\begin{gather}\label{eq:DSCC2017_5}
\Delta V(i) 
=(\rho-1)s^2(i)+\rho_{\beta}\tilde{\beta}(i)\Big(\frac{Ts(i)ax(i)}{\rho_{\beta}}+\Delta\tilde{\beta}(i)\Big)\\
+\rho_{\alpha}\tilde{\alpha}(i)\Big(\frac{Ts(i)x(i)}{\rho_{\alpha}}+\Delta\tilde{\alpha}(i)\Big) \nonumber
\end{gather}
If $\hat{\beta}$ and $\hat{\alpha}$ are updated according to the following rules: \vspace{-0.5cm}
\begin{subequations} \label{eq:DSCC2017_6}
\begin{align}
\hat{\beta}(i+1)=\hat{\beta}(i)+\frac{Ts(i)ax(i)}{\rho_{\beta}} \\
\hat{\alpha}(i+1)=\hat{\alpha}(i)+\frac{Ts(i)x(i)}{\rho_{\alpha}} 
\end{align} 
\end{subequations}
then, Eq.~(\ref{eq:DSCC2017_5}) becomes:\vspace{-0.5cm}
\begin{gather}\label{eq:DSCC2017_7}
\Delta V(i) 
=-(1-\rho)s^2(i)
\end{gather}
As can be seen, since $0<\rho<1$, Eq.~(\ref{eq:DSCC2017_7}) is negative semi-definite. This means that the positive definite Lyapunov function $V$ has a negative semi-definite difference function ($\Delta V$). Thus, according to the Lyapunov stability theorem and the new Invariance Principle for discontinuous systems~\cite{barkana2015new,Slotine,Amini_DSC,Selina_PhD}, the asymptotic stability of the closed loop controller with the adaptation laws in Eq.~(\ref{eq:DSCC2017_6}) is guaranteed. This ensures the finite time convergences of the first order sliding function ($s$) and the unknown parameter estimation errors ($\tilde{\beta},~\tilde{\alpha}$) to zero. 

The performed Lyapunov stability can be extended to higher order systems. For a second order system ($r=2$), two scalar positive definite Lyapunov functions ($V_1,~V_2$) can be defined for each of the system's states. If one can show that both Lyapunov functions have negative semi-definite difference functions, the overall stability of the second order linear system can be concluded. For a SISO first order DSMC, the sliding function of the first state ($s_1$) can be obtained as follows according to Eq.~(\ref{eq:Generic_linear_3}): \vspace{-0.5cm}
\begin{gather}\label{eq:DSCC2017_8}
s_1(i+1)=\rho_1 s_1(i)+\\
T\Big((\tilde{\beta}_{11}a_{11}+\tilde{\alpha}_{11})x_1(i)+(\tilde{\beta}_{12}a_{12}+\tilde{\alpha}_{12})x_2(i)\Big) \nonumber
\end{gather}
For $s_1$, the following positive definite Lyapunov function is introduced:\vspace{-0.6cm}
\begin{gather} \label{eq:DSCC2017_9}
V_1(i)=\frac{1}{2}{s_1^2}(i)+\frac{1}{2}\rho_{\beta_{11}}\tilde{\beta}_{11}^2(i)+\frac{1}{2}\rho_{\beta_{12}}\tilde{\beta}_{12}^2(i)+\\
\frac{1}{2}\rho_{\alpha_{11}}\tilde{\alpha}_{11}^2(i)+\frac{1}{2}\rho_{\alpha_{12}}\tilde{\alpha}_{12}^2(i) \nonumber
\end{gather}
Similar to the first order linear system stability analysis, it can be easily shown that
~the Lyapunov difference function ($\Delta V_1$) becomes:\vspace{-0.5cm}
\begin{gather}\label{eq:DSCC2017_10}
\Delta V_1(i) 
=(\rho_1-1)s_1^2(i)+\\\rho_{\beta_{11}}\tilde{\beta}_{11}(i)\Big(\frac{T.s_1(i).a_{11}.x_{1}(i)}{\rho_{\beta_{11}}}+\Delta\tilde{\beta}_{11}(i)\Big)+\nonumber\\
\rho_{\alpha_{11}}\tilde{\alpha}_{11}(i)\Big(\frac{T.s_1(i).x_1(i)}{\rho_{\alpha_{11}}}+\Delta\tilde{\alpha}_{11}(i)\Big)+\nonumber\\
\rho_{\beta_{12}}\tilde{\beta}_{12}(i)\Big(\frac{T.s_1(i).a_{12}.x_{2}(i)}{\rho_{\beta_{12}}}+\Delta\tilde{\beta}_{12}(i)\Big)+\nonumber\\
\rho_{\alpha_{12}}\tilde{\alpha}_{12}(i)\Big(\frac{T.s_1(i).x_2(i)}{\rho_{\alpha_{12}}}+\Delta\tilde{\alpha}_{12}(i)\Big) \nonumber
\end{gather}
As can be seen, if the adaptation laws from Eq.~(\ref{eq:Generic_linear_4}) are used to update $\tilde{\beta}_{11}$, $\tilde{\beta}_{12}$, $\tilde{\alpha}_{11}$, and $\tilde{\alpha}_{12}$ in Eq.~(\ref{eq:DSCC2017_10}), the Lyapunov difference function becomes $\Delta V_1(i)=-(1-\rho_1)s_1^2(i)$, which fulfills the required negative semi-definite condition for $\Delta V_1$ and guarantees the finite-time zero convergence of $s_1$, $\tilde{\beta}_{11}$, $\tilde{\beta}_{12}$, $\tilde{\alpha}_{11}$, and $\tilde{\alpha}_{12}$. The same conclusions can be reached for higher order systems ($r>2$) by utilizing the adaptation laws in Eq.~(\ref{eq:Generic_linear_4}) to update the unknown multiplicative and additive parameters within the linear system dynamics.~$\blacksquare$

Overall, the control input of an \textit{adaptive first order DSMC} with incorporated implementation imprecisions and adaptation laws from Eq.~(\ref{eq:Generic_linear_4}) becomes:\vspace{-0.5cm}
\begin{gather}\label{eq:Generic_linear_5}
^{1DSMC}\mathbf{U}^{mod}_{adaptive}(i)=
\frac{\mathbf{B}^{-1}}{T}\Big(\mathcal{P}_{r\times r}\mathbf{S}_{r\times 1}(i)+\mathbf{X}_{d}(i+1)
\\
-(T\hat{\mathbf{A}}
+\mathbf{I}_{r \times r})\mathbf{X}(i)
\Big)
-|\mu_{\mathbf{U}}|\times sat(\mathbf{S}_{r \times 1}(i)) \nonumber \vspace{-0.5cm}
\end{gather} 
where: \vspace{-0.5cm}
\begin{gather}
\hat{\mathbf{A}}=\begin{bmatrix}
\hat{\beta}_{11} a_{11}+\hat{\alpha}_{11}~~ &  \hat{\beta}_{12} a_{12}+\hat{\alpha}_{12} & \cdots &  \hat{\beta}_{1r} a_{1r}+\hat{\alpha}_{1r} \\ 
\hat{\beta}_{21} a_{21}+\hat{\alpha}_{21}~~ &  \hat{\beta}_{22} a_{22}+\hat{\alpha}_{22} & \cdots &  \hat{\beta}_{2r} a_{2r}+\hat{\alpha}_{2r}\\ 
\vdots & \vdots & \ddots & \vdots \\ 
\hat{\beta}_{r1} a_{r1}+\hat{\alpha}_{r1}~~ &  \hat{\beta}_{r2} a_{r2}+\hat{\alpha}_{r2} & \cdots &  \hat{\beta}_{rr} a_{rr}+\hat{\alpha}_{rr} \end{bmatrix} 
\end{gather}
\vspace{-0.65cm}
\subsection{Second Order Adaptive DSMC} \label{SecondOrderDSMC}
As we have shown in our previous work~\cite{Amini_ACC2017}, the concept of second order sliding mode is an effective solution to (i) minimize the high frequency oscillation due to chattering phenomena, and (ii) enhance the first order DSMC robustness characteristics against implementation imprecisions. The better performance of the second
order DSMC can be traced in driving the second order derivative (difference function) of the system states to the sliding manifold, in addition to the zero convergence of the sliding variable itself. The second order sliding mode for a continuous-time system is determined by the following equalities~\cite{Amini_ACC2017,salgado2004robust}:
\vspace{-0.4cm}
\begin{gather}\label{eq:DSCC2017_11}
\mathbf{S}(t,x)=\dot{\mathbf{S}}(t,x)=0
\end{gather}
In order to convert the second order sliding mode to a first order one, a new sliding variable ($\xi$) is defined according to $s$ and $\dot{s}$:\vspace{-0.4cm}
\begin{gather}\label{eq:DSCC2017_12}
\xi(t,x)=\dot{s}(t,x)+\lambda s(t,x),~\lambda>0
\end{gather}
$\xi$ is the sliding surface of a system with a relative order equal to one, in which the input is $\dot{u}$ and output is $\xi(t,x)$~\cite{sira1990structure}. Introduction of $\xi$ helps to follow the first order DSMC design procedure in Section~\ref{FirstOrderDSMC} in order to establish the second order DSMC formulation. According to Eq.~(\ref{eq:DSCC2017_12}), for the discrete-time system, the second order sliding function vector ($\boldsymbol{\Xi}=[\xi_1,...,\xi_r]^\intercal$) is defined as:
\vspace{-0.5cm}
\begin{gather}\label{eq:DSCC2017_13}
\boldsymbol{\Xi}(i)=\mathbf{S}(i+1)+{\Phi}\mathbf{S}(i),
\end{gather}
where $\mathbf{S}(i+1)$ is calculated with respect to Eq.~(\ref{eq:StageIII_9}) and $\Phi\in{\mathbb{R}^{r\times r}}$ is the \textit{positive definite} matrix of the second order sliding mode gains~\cite{Amini_ACC2017}. The second order DSMC control input is calculated by solving the following equalities in discrete-time~\cite{mihoub2009real}:
\vspace{-0.45cm}
\begin{gather}\label{eq:DSCC2017_14}
\boldsymbol{\Xi}(i+1)=\boldsymbol{\Xi}(i)=0
\end{gather}
Applying Eq.~(\ref{eq:DSCC2017_14}) to the linear system in Eq.~(\ref{eq:Generic_linear_2}), in the absence of model and ADC uncertainties, results in the following control input for the second order DSMC ($^{2DSMC}\mathbf{U}$):
\vspace{-0.35cm}
\begin{gather}\label{eq:DSCC2017_15}
^{2DSMC}\mathbf{U}(i)=\\
\mathbf{B}^{-1}\left(\frac{1}{T}\left[(\Phi+\mathbf{I})\mathbf{X}(i)-\Phi\mathbf{X}_d(i)+\mathbf{X}_d(i+1)\right]-\mathbf{A}\mathbf{X}(i)\right)\nonumber \vspace{-0.5cm}
\end{gather}
The control input of the baseline second order DSMC (Eq.~(\ref{eq:DSCC2017_15})) can be modified ($^{2DSMC}\mathbf{U}^{mod}$) against sampling and quantization imprecisions by inclusion of the propagated ADC uncertainties ($\mu_{\mathbf{U}}$)~\cite{Amini_ACC2017}:
\vspace{-0.5cm}
\begin{gather}\label{eq:DSCC2017_16}
^{2DSMC}\mathbf{U}^{mod}(i)=^{2DSMC}\mathbf{U}(i)-|\mu_{\mathbf{U}}(i)|\times sat({\Xi}(i-1)) \vspace{-0.5cm}
\end{gather}
where $\mu_{\mathbf{U}}$ is calculated according to Eq.~(\ref{eq:StageIII_14}). By applying Eq.~(\ref{eq:DSCC2017_16}) to the linear system (Eq.~(\ref{eq:Generic_linear_2})) in the presence of the unknown multiplicative and additive terms, the vector ($\Xi$) of the second order sliding variables ($\xi$) becomes: \vspace{-0.6cm}
\begin{equation}\label{eq:DSCC2017_16}
\begin{split}
\Xi_{r\times 1}(i)=
T\begin{bmatrix}
\tilde{\beta}_{11} a_{11}+\tilde{\alpha}_{11}~~ &\cdots & \tilde{\beta}_{1r} a_{1r}+\tilde{\alpha}_{1r} \\ 
\tilde{\beta}_{21} a_{21}+\tilde{\alpha}_{21}~~& \cdots & \tilde{\beta}_{2r} a_{2r}+\tilde{\alpha}_{2r}\\ 
\vdots & \ddots & \vdots \\ 
\tilde{\beta}_{r1} a_{r1}+\tilde{\alpha}_{r1}~~& \cdots & \tilde{\beta}_{rr} a_{rr}+\tilde{\alpha}_{rr} 
\end{bmatrix}\mathbf{X}_{r\times 1}(i)
\end{split}
\end{equation} 

It can be shown that the adaptation law for handling the modeling uncertainties for the second order DSMC has the same structure of the proposed adaptation laws for the first order DSMC (Eq.~(\ref{eq:Generic_linear_4}))~\cite{Amini_ACC2017}. The detailed proof of the latter statement for the adaptation laws of a second order DSMC design is not discussed here due to the page limit. However, similar to the first order DSMC, a brief sketch of the proof is as follows: (i) starting with a first order system ($r=1$), the following Lyapunov function is used to begin with the closed loop system stability analysis:
\vspace{-0.7cm}
\begin{gather}\label{eq:D2SMC_8}
V(k)=\frac{1}{2}\Big({s}^2(i+1)+\varphi{s}^2(i)\Big)\\
+\frac{1}{2}\rho_{\beta}\Big(\tilde{\beta}^2(i+1)+\varphi\tilde{\beta}^2(i)\Big)\nonumber\\
+\frac{1}{2}\rho_{\alpha}\Big(\tilde{\alpha}^2(i+1)+\varphi\tilde{\alpha}^2(i)\Big),\nonumber
\end{gather}
(ii) by calculation of the Lyapunov difference function ($\Delta V$) according to Eq.~(\ref{eq:D2SMC_8}), it can be shown that by utilizing Eq.~(\ref{eq:Generic_linear_4}) to estimate the unknown parameters, and choosing $0<\varphi<1$, the asymptotic stability of the closed-loop system is guaranteed, and (iii) finally, the Lyapunov stability analysis can be extended to higher order SISO/MIMO systems ($r>1$), and it can be shown that in addition to the rules in Eq.~(\ref{eq:Generic_linear_4}) as the adaptation laws, the matrix of second order DSMC gain ($\Phi$) should be $\textit{symmetric positive}$ to ensure the stability of the system. Similar to the first order DSMC, when $\Phi$ is diagonal, the second order DSMC is SISO, and when the off-diagonal elements of $\Phi$ are not zero, there is coupling between inputs and outputs of the closed loop system, and the controller has a MIMO structure. \vspace{-0.65cm}
\section{Case Study: DC Motor Speed Control} \label{CS_DCMotor}
\begin{figure}[h!]
\begin{center}
\includegraphics[angle=0,width= \columnwidth]{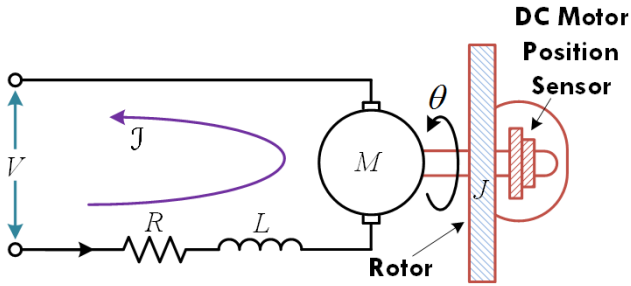} \vspace{-0.6cm}
\caption{\label{fig:DCMotor_Schematic}Schematic of the modeled DC motor.} \vspace{-1cm}
\end{center}
\end{figure}

DC motors are common actuators for control applications that require rotary and transitional motions. The electric equivalent circuit of the armature and the free-body diagram of the rotor are shown in Figure~\ref{fig:DCMotor_Schematic}. For speed regulation of a DC motor, the control input is voltage ($V$) to the motor's armature and the output is rotation speed (${\theta}$) of the shaft. Assuming a constant magnetic field and linear relationship between motor torque and armature current ($\mathcal{I}$), and by choosing the rotor speed and current as the state variables, the following linear time-invariant state-space representation can be used to describe the dynamics of the DC motor~\cite{DCMotor_ref}: \vspace{-0.6cm}
%
\begin{subequations} \label{eq:DC_motor_linear_DSMC}
\begin{align}
\theta(i+1)=T\left(\frac{k_m}{J}{\mathcal{I}}(i)-\frac{k_f}{J}\theta(i)+\frac{1}{J} \Gamma\right)+\theta(i) \\
{\mathcal{I}}(i+1)=T\left(-\frac{k_b}{L}\theta(i)-\frac{R}{L}{\mathcal{I}}(i)+\frac{1}{L}V(i)\right)+{\mathcal{I}}(i)
\end{align}
\end{subequations}
where $J$ is the rotor's moment of inertia, $R$ is the electrical resistance, $L$ is the electric inductance, $\Gamma$ is the torque on the rotor, $k_f$ is the mechanical damping, $k_m$ is the motor torque constant, and $k_b$ is the electromotive force constant. The DC motor model constants are listed in 
the Appendix. \vspace{-0.85cm}

\begin{figure}[h!]
\begin{center}
\includegraphics[angle=0,width= \columnwidth]{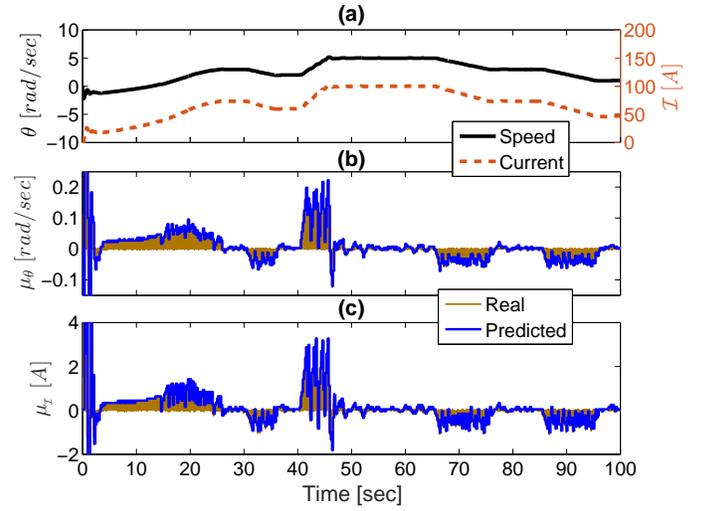} \vspace{-0.65cm}
\caption{\label{fig:DCMotor_SpeedU_PU}Comparison of actual and predicted uncertainties due to sampling and quantization. (a) measured signals (i.e., signal after ADC), (b) actual and predicted uncertainties on measured speed signal, and (c) actual and predicted uncertainties on current signal. Test conditions: 200~$ms$ sampling time and 10-$bit$ quantization level.}\vspace{-0.75cm}
\end{center}
\end{figure}

Performance of the employed uncertainty prediction technique is investigated for the DC motor, under 200 $ms$ of sampling time and a 10-$bit$ of quantization level. $\theta$ and $\mathcal{I}$ are feedback signals that go through ADC before going to the DC motor controller. Simulations are done in MATLAB Simulink$^{\textregistered}$~that allows for testing the controller in a model-in-the-loop (MIL) platform against sampling and quantization imprecisions. Figure~\ref{fig:DCMotor_SpeedU_PU} shows the ADC uncertainty prediction results on the shaft speed and the current of the armature circuit. In this figure, the term ``real" denotes the signal after ADC in the MIL setup in MATLAB. The accuracy of the uncertainty prediction technique is shown in Table~\ref{tab:DC_motor_table_StageIII_2} in terms of mean error and standard deviation of the errors. The small error values in Table~\ref{tab:DC_motor_table_StageIII_2} confirm the capability of the proposed method to estimate ADC uncertainty on measured signals for the DC motor case study. In the next section, the predicted uncertainties on speed and current are employed to calculate propagated uncertainties on control inputs which will be directly utilized to design a DSMC for the DC motor speed regulation under implementation imprecisions.

{\small
\linespread{0.8}
\begin{table}[h!]
\begin{center}
\caption{ Mean ($\bar{e}$) and standard deviation ($\sigma_e$) of ADC uncertainty prediction errors for the DC motor case study.}
\label{tab:DC_motor_table_StageIII_2} 
    \begin{tabular}{lcc}
       \hline\hline
        & {$\bar{e}$} & {$\sigma_e$} \\
       \hline
       $Motor~Speed,~\theta~[rad/sec]$ &  0.004 & 0.094  \\ \vspace{0.05cm}
       $Current,~{\mathcal{I}}~[A]$ 	 &  0.073 & 0.823  \\
        \hline\hline 
\end{tabular}
\end{center}\vspace{-0.65cm}
\end{table}
}
%

%

In the absence of model uncertainties, the SISO first and second order DSMCs are designed to regulate the DC motor rotational speed with respect to the desired speed profile ($\theta_d$) under implementation imprecisions. The first sliding surface is defined as the error in tracking the desired speed profile ($s_{1}=\theta-\theta_d$). Since there is no direct control input on DC motor rotational speed, $\mathcal{I}_d$ is defined as the synthetic control input for controlling the shaft speed. $\mathcal{I}_d$ is used to define the second sliding surface ($s_{2}=\mathcal{I}-\mathcal{I}_d$) in which the control input is voltage. Upon incorporating the predicted ADC uncertainties on control signals into the DSMC design, synthetic ($\mathcal{I}_d$) and physical ($V$) control inputs of the first order DSMC are calculated according to Eq.~(\ref{eq:StageIII_21}) and Eq.~\eqref{eq:DC_motor_linear_DSMC}: \vspace{-0.5cm}
%
%
\begin{equation}\label{eq:DC_motor_DSMC_modified}
 \begin{split}
 ^{1DSMC}{\mathcal{I}}^{mod}_d(i)=\frac{J}{k_m}\left(\frac{1}{T}\left[\rho_1(\theta(i)-\theta_d(i))+\theta_d(i+1)\right.\right.~~~~~~\\
 \left.\left.-\theta(i)\right]+\frac{k_f}{J}\theta(i)-\frac{1}{J}\Gamma\right)-|\mu_{{\mathcal{I}}_d}(i)|sat(\theta(i)-\theta_d(i)) 
 \end{split}
 \end{equation} 
\vspace{-1.1cm}

 \begin{equation}\label{eq:DC_motor_DSMC_modified2}
 \begin{split}
^{1DSMC}V^{mod}(i)=L\left(\frac{1}{T_s}\left[\rho_2({\mathcal{I}}(i)-{\mathcal{I}}_d(i))+{\mathcal{I}}_d(i+1) \right.\right.\\ 
 \left.\left.-{\mathcal{I}}(i)\right]+\frac{k_b}{L}\theta(i)+\frac{R}{L}{\mathcal{I}}(i)\right)-|\mu_{_V}(i)|sat({\mathcal{I}}(i)-{\mathcal{I}}_d(i))
 \end{split}
 \end{equation} 
where $\mu_{{\mathcal{I}}_d}$ and $\mu_{_V}$ denote the estimations of the propagated ADC uncertainties on control signals (Eq.~(\ref{eq:StageIII_15_1})):\vspace{-0.6cm}
\begin{gather}\label{eq:DC_motor_DSMC_mu}
\mu_{{\mathcal{I}}_d}(i)=\frac{J}{Tk_m}\Big((\rho_1-1)\mu_{\theta}(i)\Big)+\frac{k_f}{k_m}\mu_{\theta}(i) \\
\mu_{_V}(i)=\frac{L}{T}\Big((\rho_2-1)\mu_{\mathcal{I}}(i)\Big)+k_b\mu_\theta(i)+R\mu_{\mathcal{I}}(i) 
\end{gather} 
$\mu_\theta$ and $\mu_{\mathcal{I}}$ are predicted uncertainties on measured signals that are calculated using Eq.~(\ref{eq:StageIII_7}). These uncertainties were previously shown in Figure~\ref{fig:DCMotor_SpeedU_PU}. The state equations of the DC motor model with unknown multiplicative and additive terms becomes: \vspace{-0.4cm}
\begin{gather}\label{eq:DC_motor_linear_DSMC_Uncertain}
\theta(i+1)=T\left([-{\beta}_{11}\frac{k_f}{J}+\alpha_{11}]\theta(i)+[\alpha_{12}]{\mathcal{I}}(i)\right)\\
+T\left(\frac{k_m}{J}{\mathcal{I}}v(i)+\frac{1}{J} \Gamma\right) +\theta(i)\nonumber \\
{\mathcal{I}}(i+1)=T\left([-\beta_{21}\frac{k_b}{L}+\alpha_{21}]\theta(i)+[-\beta_{22}\frac{R}{L}+\alpha_{22}]{\mathcal{I}}(i)\right)\\
+\frac{T}{L}V(i)+{\mathcal{I}}(i) \nonumber   
\end{gather}

Four additive ($\alpha_{pq}$) and three multiplicative ($\beta_{pq}$) unknown parameters are derived to their nominal values by solving the adaptation laws in Eq.~(\ref{eq:Generic_linear_4}). The final adaptive first order SISO DSMC with predicted ADC uncertainties yields: \vspace{-0.5cm}
\begin{equation}\label{eq:DC_motor_adaptiveDSMC_modified}
 \begin{split}
^{1DSMC}{\mathcal{I}}^{mod}_{d,adaptive}(i)=\frac{J}{k_m}\left(\frac{1}{T}\left[\rho_1(s_1(i))+\theta_d(i+1)-\theta(i)\right]+\right.\\
 \left.[\hat{\beta}_{11}\frac{k_f}{J}+\hat{\alpha}_{11}]\theta(i)+[\hat{\alpha}_{21}]{\mathcal{I}}(i)-\frac{1}{J}\Gamma\right)-|\mu_{{\mathcal{I}}_d}(i)|sat(s_1(i)) 
 \end{split}
 \end{equation} 
 \vspace{-1.5cm}
 \begin{equation}\label{eq:DC_motor_adaptiveDSMC_modified2}
 \begin{split}
^{1DSMC}V^{mod}_{adaptive}(i)=L\left(\frac{1}{T}\left[\rho_2(s_2(i))+{\mathcal{I}}_d(i+1)-{\mathcal{I}}(i)\right] \right.\\ 
 \left.+[\hat{\beta}_{21}\frac{k_b}{L}+\hat{\alpha}_{21}]\theta(i)+[\hat{\beta}_{22}\frac{R}{L}+\hat{\alpha}_{22}]{\mathcal{I}}(i)\right)-|\mu_{_V}(i)|sat(s_2(i))
 \end{split}
 \end{equation} 

In a similar manner to the first order DSMC design for the DC motor case study, the control inputs of the adaptive second order SISO DSMC with predicted ADC uncertainties are as follows:\vspace{-0.5cm}
\begin{equation}\label{eq:DC_motor_adaptive2DSMC_modified}
 \begin{split}
^{2DSMC}{\mathcal{I}}^{mod}_{d,adaptive}(i)=\frac{J}{k_m}\left(\frac{1}{T}\left[-\varphi_1(s_1(i))+\theta_d(i+1)-\theta(i)\right]\right.\\
 \left.+[\hat{\beta}_{11}\frac{k_f}{J}+\hat{\alpha}_{11}]\theta(i)+[\hat{\alpha}_{21}]{\mathcal{I}}(i)-\frac{1}{J}\Gamma\right)-|\mu_{{\mathcal{I}}_d}(i)|sat(\xi_1(i-1)) 
 \end{split}
 \end{equation} 
 \vspace{-1.5cm}
 \begin{equation}\label{eq:DC_motor_adaptiveDSMC_modified2}
 \begin{split}
^{2DSMC}V^{mod}_{adaptive}(i)=L\left(\frac{1}{T}\left[-\varphi_2(s_2(i))+{\mathcal{I}}_d(i+1)-{\mathcal{I}}(i)\right] \right.\\ 
 \left.+[\hat{\beta}_{21}\frac{k_b}{L}+\hat{\alpha}_{21}]\theta(i)+[\hat{\beta}_{22}\frac{R}{L}+\hat{\alpha}_{22}]{\mathcal{I}}(i)\right)-|\mu_{_V}(i)|sat(\xi_2(i))
 \end{split}
 \end{equation} \vspace{-0.75cm}

\begin{figure}[h!]
\begin{center}
\includegraphics[angle=0,width= \columnwidth]{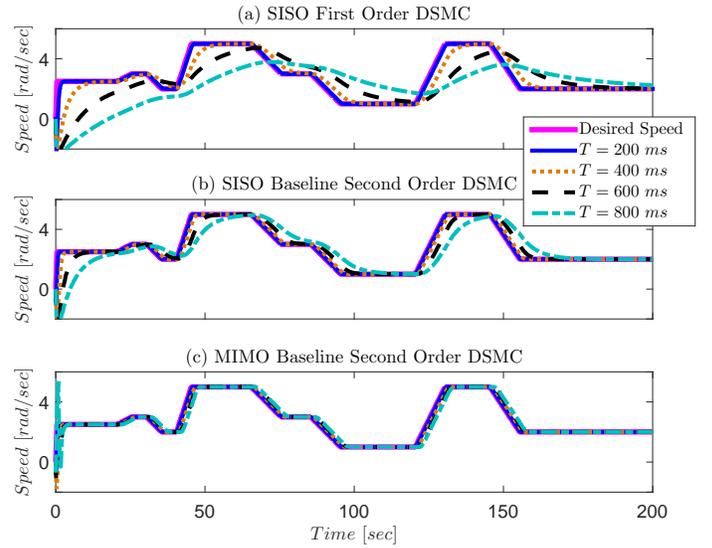} \vspace{-0.75cm}
\caption{\label{fig:DC_2DSMC_SamplingComparison}Comparison among the speed tracking results of the first and second order DSMCs for different sampling rates and quantization level of 16-$bit$: (a) SISO first order DSMC, (b) SISO second order DSMC, (c) MIMO second order DSMC. No model uncertainty is applied.} \vspace{-0.85cm}
\end{center}
\end{figure}

The SISO second order DSMC from Eq.~(\ref{eq:DC_motor_adaptive2DSMC_modified}) and (\ref{eq:DC_motor_adaptiveDSMC_modified2}) can be converted into a MIMO structure via the second order sliding mode tuning gains ($\varphi$), in which the off-diagonal elements of $\Phi$ are chosen to be non-zero to reflect the coupling between the states of the DC motor model in the controller structure. In the absence of model uncertainties ($\beta_{pq}=1,~\alpha_{pq}=0$), Figure~\ref{fig:DC_2DSMC_SamplingComparison} shows the comparison between the SISO first, and SISO/MIMO second order DSMCs for different sampling rates. Shannon's sampling theorem criteria states that the sampling frequency must be at least twice the maximum frequency of the measured analog signal. As long as Shannon's sampling theorem is satisfied, increasing the sampling time helps to reduce the computation cost. Although at lower sampling rates (e.g., $200~ms$) all the controllers show similar performances, by increasing the sampling rate to $800~ms$, the higher robustness characteristics of both SISO and MIMO second order DSMCs in comparison with the first order controller is revealed. The comparison results show that the SISO second order DSMC is improving the tracking errors by 69\% on average for different sampling rates, compared to the first order controller. Moreover, {except for the overshoot at the beginning which becomes larger as the sampling time increases}, it can be seen that the MIMO second order controller is barely affected by the sampling time increase, and this illustrates its strong robustness against ADC uncertainties.\vspace{-0.65 cm}

\begin{figure}[h!]
\begin{center}
\includegraphics[angle=0,width= \columnwidth]{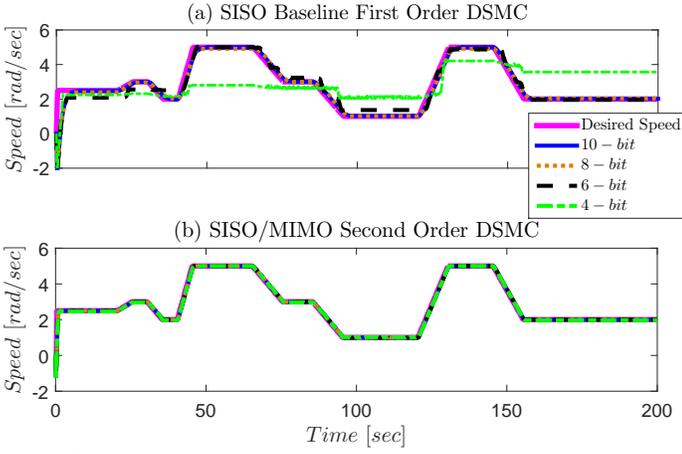} \vspace{-0.75cm}
\caption{\label{fig:DC_2DSMC_QuantizationComparison}Comparison between the first and second order DSMCs in tracking the desired speed trajectory of the DC motor for different quantization levels. (a) SISO first order DSMC, (b) SISO/MIMO second order DSMC ($T$=200$~ms$, no model uncertainty is applied).} \vspace{-0.75cm}
\end{center}
\end{figure}

Figure~\ref{fig:DC_2DSMC_QuantizationComparison} shows how the DSMC tracking performance could be affected by the quantization level of the ADC unit. As can be seen from Figure~\ref{fig:DC_2DSMC_QuantizationComparison}-a, the SISO first order DSMC is more sensitive to quantization level compared to the SISO/MIMO second order DSMC. On the other side, any changes in the quantization level from 10-$bit$ to 4-$bit$ have no effect on both SISO and MIMO second order DSMCs tracking performances.

\begin{figure}[h!]
\begin{center}
\includegraphics[angle=0,width= \columnwidth]{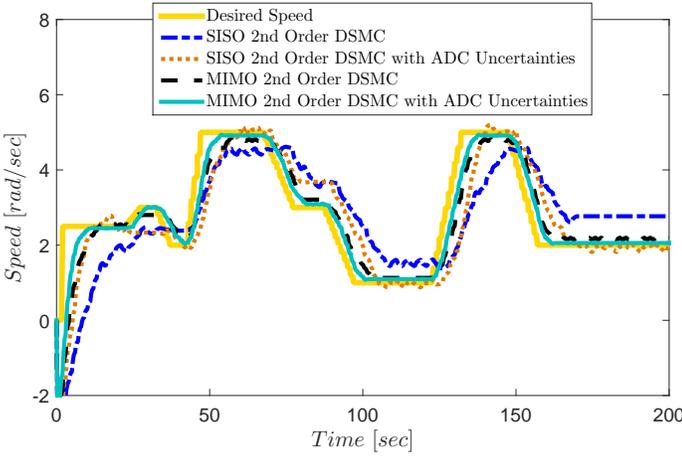} \vspace{-0.75cm}
\caption{\label{fig:DC_2DSMC_MIMOSISO}Comparison between the MIMO and SISO second order DSMCs, with and without predicted ADC uncertainties ($\mu_u$), in tracking the desired speed trajectory under extreme ADC uncertainties ($T$=1000$~ms$, quantization level=4-$bit$, no model uncertainty is applied).} \vspace{-0.85cm}
\end{center}
\end{figure}

In order to show the effectiveness of incorporating the predicted ADC uncertainties in improving the DSMC against uncertainties, the second order DSMC is evaluated under extreme sampling and quantization levels which cause high level of imprecisions at the controller I/O. By looking into Figure~\ref{fig:DC_2DSMC_MIMOSISO}, one can conclude the better tracking performance of the second order DSMC with predicted ADC uncertainties for both SISO and MIMO cases, compared to the two other controllers. 

\begin{figure}[h!]
\begin{center}
\includegraphics[angle=0,width= \columnwidth]{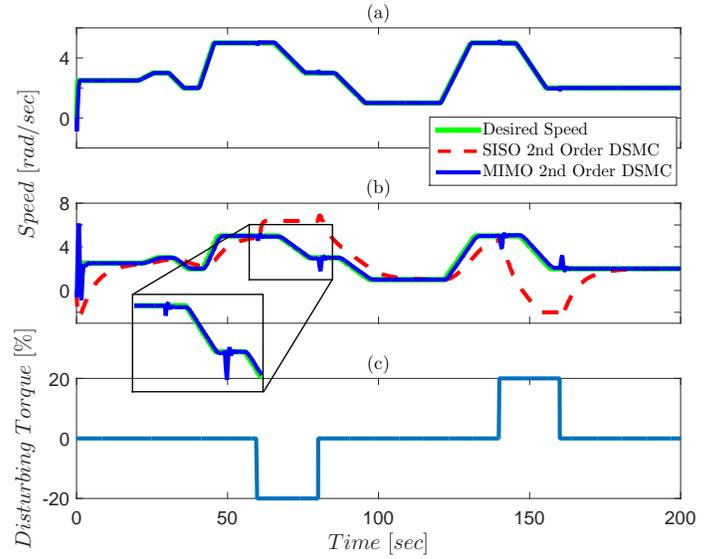} \vspace{-0.75cm}
\caption{\label{fig:DC_2DSMC_DisturbanceRejection}Disturbance rejection performance of the SISO and MIMO second order DSMCs for up to 20$\%$ sudden changes on the nominal external torque for (a) $T$=200~$ms$, quantization level=10-$bit$, (b) $T$=800$~ms$, quantization level=10-$bit$. Disturbing torque percentage is plotted in (c). No model uncertainty is applied.} \vspace{-0.8cm}
\end{center}
\end{figure}

Another interesting feature of the MIMO second order DSMC in comparison with the SISO controller is its disturbance rejection characteristics. As shown in Eq.~(\ref{eq:DC_motor_linear_DSMC_Uncertain}), it is assumed that there is a constant torque ($\Gamma$) on the shaft of the DC motor. If there will be any unknown disturbing torque on the motor shaft, the controller should reject the disturbance fast, to make sure the desired speed tracking is not affected. Figure~\ref{fig:DC_2DSMC_DisturbanceRejection} shows the disturbance rejection characteristics of the SISO and MIMO second order DSMCs, under different sampling and quantization levels. The disturbing torque, shown in Figure~\ref{fig:DC_2DSMC_DisturbanceRejection}-c, is defined with respect to the constant nominal torque on the DC motor shaft, and it generates up to 20\% disturbing torque load. Figure~\ref{fig:DC_2DSMC_DisturbanceRejection}-a shows that when the sampling rate is fast (200 $ms$), both SISO and MIMO second order DSMCs reject the disturbing torque effect very quickly. However, upon increasing the sampling rate from 200~$ms$ to 800~$ms$ (Figure~\ref{fig:DC_2DSMC_DisturbanceRejection}-b), while the MIMO controller still shows accurate speed tracking and fast disturbance rejection results, the SISO controller fails to reject the disturbing torque impact and the tracking performance is affected significantly. 

\begin{figure}[h!]
\begin{center}
\includegraphics[angle=0,width= \columnwidth]{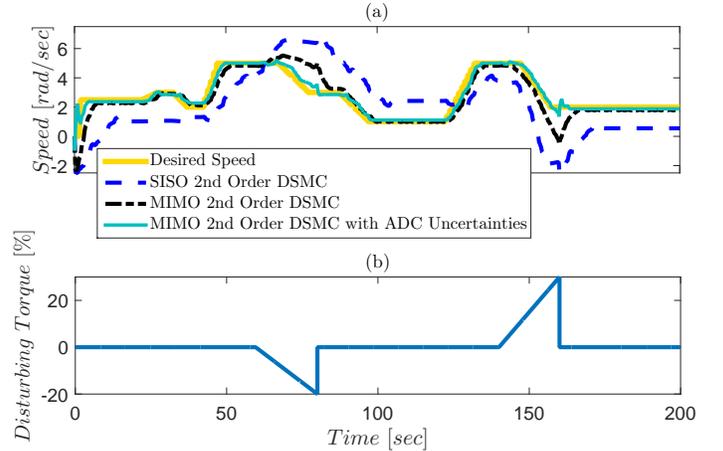} \vspace{-0.75cm}
\caption{\label{fig:DC_2DSMC_DisturbanceRejection_withADC}Disturbance rejection performance of the SISO and MIMO second order DSMCs, with and without predicted ADC uncertainties ($\mu_u$), for up to 30$\%$ changes on the nominal external torque under extreme ADC uncertainties of $T$=1000$~ms$ and quantization level=4-$bit$, and (b) disturbing torque percentage. No model uncertainty is applied.} \vspace{-0.85cm}
\end{center}
\end{figure}

Figure~\ref{fig:DC_2DSMC_DisturbanceRejection_withADC}-a shows that at the extreme sampling rate (1000~$ms$) and quantization levels (4-$bit$), even the MIMO second order baseline DSMC fails to reject the disturbing torque (Figure~\ref{fig:DC_2DSMC_DisturbanceRejection_withADC}-b). It is not a surprise that the disturbance rejection characteristics of the controller is weakened when the external uncertainties become larger. Figure~\ref{fig:DC_2DSMC_DisturbanceRejection_withADC}-b shows that incorporation of the predicted ADC uncertainties (inclusion of the switching control input) helps to improve both tracking performance and disturbing torque rejection results for the MIMO second order DSMC, compared to the baseline controllers. \vspace{-0.65cm}

\begin{figure}[h!]
\begin{center}
\includegraphics[angle=0,width=\columnwidth]{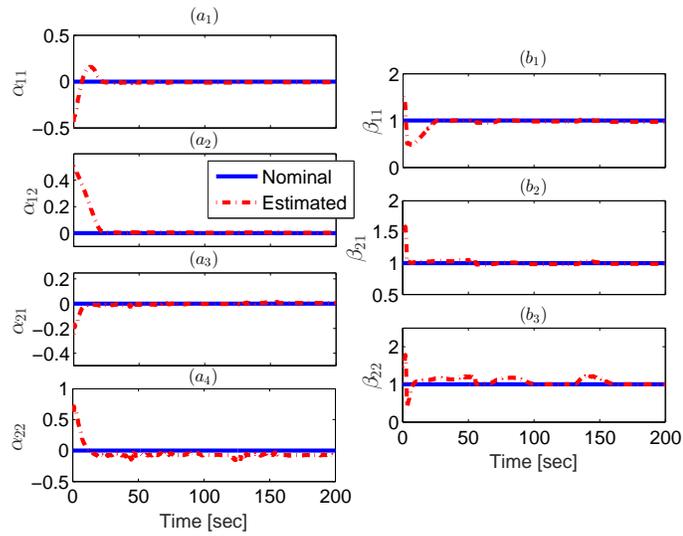}\vspace{-0.3cm}
\caption{\label{fig:DC_combined_param} Convergence results of the unknown additive ($\alpha$) and multiplicative ($\beta$) terms in the DC motor model~(200~$ms$ sampling time and 10-$bit$ quantization level).} \vspace{-0.8cm}
\end{center}
\end{figure}

In the next step, the multiplicative and additive uncertainty terms ($\beta,~\alpha$) are included into the DC motor model, by which 50\% uncertainty is introduced on the plant's parameters via each of the unknown terms. Figure~\ref{fig:DC_combined_param} shows the estimation results for four additive ($\alpha_{pq}$) and three multiplicative ($\beta_{pq}$) unknown parameters, which are steered towards their nominal values by solving the adaptation laws in Eq.~(\ref{eq:Generic_linear_4}). As can be seen, the proposed adaptation laws are able to remove the uncertainties within the model by more than 90\%. In section~\ref{sec:Results}, the final designed first and second order DSMCs are tested in real-time  by using a processor-in-the-loop (PIL) test setup on an actual ECU.\vspace{-0.82cm}

\vspace{-0.95cm}
\section{DSMC Real-Time Verification} \label{sec:Results}
\begin{figure}[h!]
\begin{center}
\includegraphics[angle=0,width= \columnwidth]{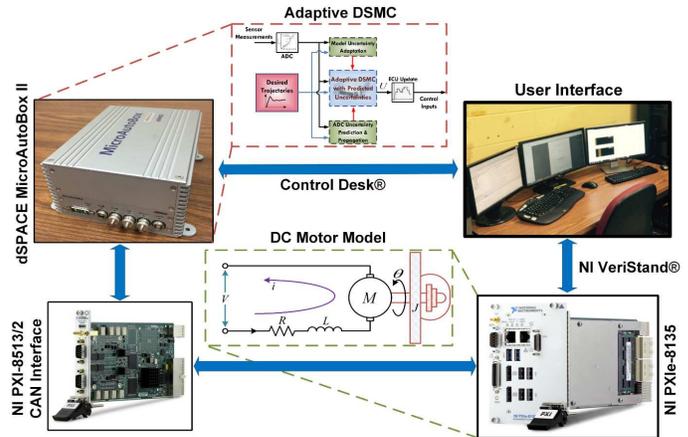} \vspace{-0.75cm}
\caption{\label{fig:HIL_PXI_Schematic}Schematic of the PIL setup for real-time DSMC verification.} \vspace{-0.6cm}
\end{center}
\end{figure}

In this section, the performance of the designed first and second order adaptive DSMCs for the DC motor speed tracking problem is evaluated on an actual ECU within a PIL setup shown in Figure~\ref{fig:HIL_PXI_Schematic}. As can be seen from Figure~\ref{fig:HIL_PXI_Schematic}, the PIL setup is equipped with two processors: (i) National Instrument (NI) PXI processor (NI PXIe-8135), and (ii) dSPACE MicroAutoboxII (MABX). The PXI processor is used for implementing the model of the DC motor, while the algorithms of the adaptive DSMCs are implemented inside the MABX. The MABX represents the actual ECU, in which the signals at the controller I/Os are sampled and quantized at a desired level by using embedded sampler and quantizer blocks. The generated C-code of the adaptive first and second order DSMCs along with the adaptation and uncertainty prediction mechanisms are implemented into MABX via dSPACE Control Desk$^{\textregistered}$~software. Real-time test configuration is conducted using NI VeriStand$^{\textregistered}$~software on an interface desktop computer. The desired speed trajectory tracking performance of the controllers are studied in real-time under $T$=200~$ms$ and quantization level=10~$bit$ on the feedback and control signals, in the presence of additive and multiplicative types of model uncertainties inside the model-based controller structure. 

\begin{figure}[h!]
\begin{center}
\includegraphics[angle=0,width= \columnwidth]{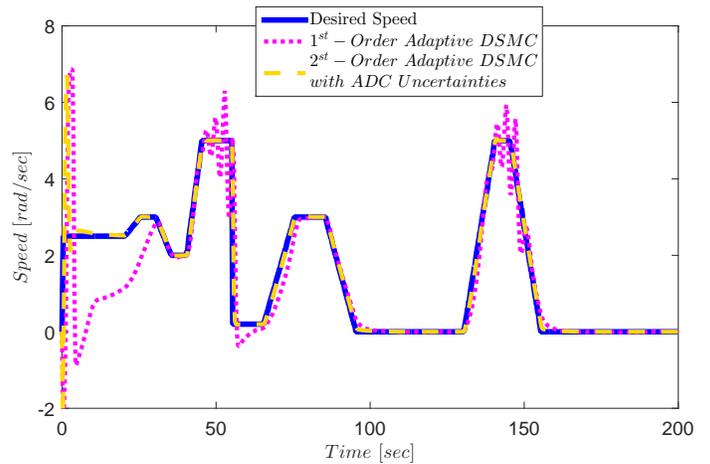} \vspace{-0.75cm}
\caption{\label{fig:HILTesting1} Real-time verification result of the DC motor speed control using the adaptive first order SISO DSMC and the modified adaptive second order MIMO DSMC with predicted ADC uncertainties. Test conditions: sampling time: 200~$ms$, quantization level: 10-$bit$, 50\% model uncertainty is applied.} \vspace{-0.8cm}
\end{center}
\end{figure}

The performance of the adaptive SISO first order DSMC, and adaptive MIMO second order DSMC with incorporated implementation uncertainties are shown in Figure~\ref{fig:HILTesting1} for tracking the desired speed profile. {Although the adaptation law for both first and second order DSMCs are the same, the adaptation mechanism for the first order DSMC is affected by the sampling imprecisions considerablly, e.g., large oscillations from 45 to 50 seconds. On the other side, the adaptive second order MIMO DSMC with predicted ADC uncertainties is able to significantly improve the tracking performance by 60\% compared to the first order adaptive SISO DSMC.}   

\vspace{-0.6cm}
\section{Summary and Conclusions}   \label{sec:Conclusion}
A new adaptive SISO/MIMO formulation of the first and second order DSMCs was presented for a general class of uncertain linear systems under implementation imprecisions and modeling uncertainties. The new formulation includes the knowledge of the ADC uncertainties via an online sampling and quantization imprecisions prediction mechanism, and adaptation laws to mitigate the modeling uncertainties via a Lyapunov stability argument. The proposed first and second order DSMCs were studied for a DC motor speed tracking control problem. The simulation and real-time PIL testing tracking results showed that: \vspace{-0.2cm}
\begin{itemize}
\item[$\bullet$] In the absence of the modeling uncertainty,  up to 84\% improvements can be achieved by using the MIMO second order DSMC compared to the SISO first order DSMC, under ADC uncertainties.
\item[$\bullet$] In the presence of the modeling uncertainty, it was observed that the adaptation mechanism is able to remove the errors in the DC motor model by up to 90\%.
\item[$\bullet$] The adaptive second order MIMO DSMC with ADC uncertainties was able to improve the adaptive SISO first order DSMC tracking performance by up to 60\%. This improvement is achieved by (i) utilizing the MIMO structure for the controller which takes into account the dynamics coupling in the controller actions, (ii) converting the first order DSMC into a second order controller, by which $\mathbf{S}$ and $\dot{\mathbf{S}}$ are both driven to zero, and (iii) incorporation of the predicted ADC uncertainties into the controller structure which enhances the overall robustness characteristics against data sampling and quantization imprecisions.
\end{itemize}
\vspace{-0.95cm}
%
%
\begin{acknowledgment}
This material is based upon the work supported by the National Science Foundation under Grant No. 1434273. Prof. J. K. Hedrick from University of California, Berkeley, and Dr. K. Butts from Toyota Motor company of North America are gratefully acknowledged for their technical comments during the course of this study.
\end{acknowledgment} \vspace{-1cm}
\bibliographystyle{IEEEtran}
\bibliography{DSMC_UP_bib}
\vspace{-0.75cm}
\section*{Appendix: Parameters~of~the~DC~Motor~Plant~Model} \vspace{0.075cm}
$J$=0.02~[m$^2$kg], $R$=2~[$\Omega$], $L$=0.5~[H], $k_m$=0.015~[N.m/A], $k_f$=0.02~[N.m.s], and $k_b$=0.015~[V.s/rad].

\end{document}